# Characterizations of joint distributions, copulas, information, dependence and decoupling, with applications to time series

**Victor H. de la Peña**[1,*], **Rustam Ibragimov**[2,†] and **Shaturgun Sharakhmetov**[3]

*Columbia University, Harvard University and Tashkent State Economics University*

**Abstract:** In this paper, we obtain general representations for the joint distributions and copulas of arbitrary dependent random variables absolutely continuous with respect to the product of given one-dimensional marginal distributions. The characterizations obtained in the paper represent joint distributions of dependent random variables and their copulas as sums of $U$-statistics in independent random variables. We show that similar results also hold for expectations of arbitrary statistics in dependent random variables. As a corollary of the results, we obtain new representations for multivariate divergence measures as well as complete characterizations of important classes of dependent random variables that give, in particular, methods for constructing new copulas and modeling different dependence structures.

The results obtained in the paper provide a device for reducing the analysis of convergence in distribution of a sum of a double array of dependent random variables to the study of weak convergence for a double array of their independent copies. Weak convergence in the dependent case is implied by similar asymptotic results under independence together with convergence to zero of one of a series of dependence measures including the multivariate extension of Pearson's correlation, the relative entropy or other multivariate divergence measures. A closely related result involves conditions for convergence in distribution of $m$-dimensional statistics $h(X_t, X_{t+1}, \ldots, X_{t+m-1})$ of time series $\{X_t\}$ in terms of weak convergence of $h(\xi_t, \xi_{t+1}, \ldots, \xi_{t+m-1})$, where $\{\xi_t\}$ is a sequence of independent copies of $X_t'$s, and convergence to zero of measures of intertemporal dependence in $\{X_t\}$. The tools used include new sharp estimates for the distance between the distribution function of an arbitrary statistic in dependent random variables and the distribution function of the statistic in independent copies of the random variables in terms of the measures of dependence of the random variables. Furthermore, we obtain new sharp complete decoupling moment and probability inequalities for dependent random variables in terms of their dependence characteristics.

*Supported in part by NSF grants DMS/99/72237, DMS/02/05791, and DMS/05/05949.

†Supported in part by a Yale University Graduate Fellowship; the Cowles Foundation Prize; and a Carl Arvid Anderson Prize Fellowship in Economics.

[1]Department of Statistics, Columbia University, Mail Code 4690, 1255 Amsterdam Avenue, New York, NY 10027, e-mail: vp@stat.columbia.edu

[2]Department of Economics, Harvard University, 1805 Cambridge St., Cambridge, MA 02138, e-mail: ribragim@fas.harvard.edu

[3]Department of Probability Theory, Tashkent State Economics University, ul. Uzbekistanskaya, 49, Tashkent, 700063, Uzbekistan, e-mail: tim001@tseu.silk.org

## 1. Introduction

In recent years, a number of studies in statistics, economics, finance and risk management have focused on dependence measuring and modeling and testing for serial dependence in time series. It was observed in several studies that the use of the most widely applied dependence measure, the correlation, is problematic in many setups. For example, Boyer, Gibson and Loretan [9] reported that correlations can provide little information about the underlying dependence structure in the cases of asymmetric dependence. Naturally (see, e.g., Blyth [7] and Shaw [71]), the linear correlation fails to capture nonlinear dependencies in data on risk factors. Embrechts, McNeil and Straumann [22] presented a rigorous study concerning the problems related to the use of correlation as measure of dependence in risk management and finance. As discussed in [22] (see also Hu [32]), one of the cases when the use of correlation as measure of dependence becomes problematic is the departure from multivariate normal and, more generally, elliptic distributions. As reported by Shaw [71], Ang and Chen [4] and Longin and Solnik [54], the departure from Gaussianity and elliptical distributions occurs in real world risks and financial market data. Some of the other problems with using correlation is that it is a bivariate measure of dependence and even using its time varying versions, at best, leads to only capturing the pairwise dependence in data sets, failing to measure more complicated dependence structures. In fact, the same applies to other bivariate measures of dependence such as the bivariate Pearson coefficient, Kullback-Leibler and Shannon mutual information, or Kendall's tau. Also, the correlation is defined only in the case of data with finite second moments and its reliable estimation is problematic in the case of infinite higher moments. However, as reported in a number of studies (see, e.g., the discussion in Loretan and Phillips [55], Cont [11] and Ibragimov [33, 34] and references therein), many financial and commodity market data sets exhibit heavy-tailed behavior with higher moments failing to exist and even variances being infinite for certain time series in finance and economics. A number of frameworks have been proposed to model heavy-tailedness phenomena, including stable distributions and their truncated versions, Pareto distributions, multivariate $t$-distributions, mixtures of normals, power exponential distributions, ARCH processes, mixed diffusion jump processes, variance gamma and normal inverse Gamma distributions (see [11, 33, 34] and references therein), with several recent studies suggesting modeling a number of financial time series using distributions with "semiheavy tails" having an exponential decline (e.g., Barndorff–Nielsen and Shephard [5] and references therein). The debate concerning the values of the tail indices for different heavy-tailed financial data and on appropriateness of their modeling based on certain above distributions is, however, still under way in empirical literature. In particular, as discussed in [33, 34], a number of studies continue to find tail parameters less than two in different financial data sets and also argue that stable distributions are appropriate for their modeling.

Several approaches have been proposed recently to deal with the above problems. For example, Joe [42, 43] proposed multivariate extensions of Pearson's coefficient and the Kullback–Leibler and Shannon mutual information. A number of papers have focused on statistical and econometric applications of mutual information and other dependence measures and concepts (see, among others, Lehmann [52], Golan [26], Golan and Perloff [27], Massoumi and Racine [57], Miller and Liu [58], Soofi and Retzer [73] and Ullah [76] and references therein). Several recent papers in econometrics (e.g., Robinson [66], Granger and Lin [29] and Hong and White [31]) considered problems of estimating entropy measures of serial dependence in time series. In a study of multifractals and generalizations of Boltzmann-Gibbs statistics, Tsallis [75] proposed a class of generalized entropy measures that include, as a particular case, the Hellinger distance and the mutual information measure. The latter measures were used by Fernandes and Flôres [24] in testing for conditional independence and noncausality. Another approach, which is also becoming more and more popular in econometrics and dependence modeling in finance and risk management is the one based on copulas. Copulas are functions that allow one, by a celebrated theorem due to Sklar [72], to represent a joint distribution of random variables (r.v.'s) as a function of marginal distributions (see Section 3 for the formulation of the theorem). Copulas, therefore, capture all the dependence properties of the data generating process. In recent years, copulas and related concepts in dependence modeling and measuring have been applied to a wide range of problems in economics, finance and risk management (e.g., Taylor [74], Fackler [23], Frees, Carriere and Valdez [25], Klugman and Parsa [46], Patton [61, 62], Richardson, Klose and Gray [65], Embrechts, Lindskog and McNeil [21], Hu [32], Reiss and Thomas [64], Granger, Teräsvirta and Patton [30] and Miller and Liu [58]). Patton [61] studied modeling time-varying dependence in financial markets using the concept of conditional copula. Patton [62] applied copulas to model asymmetric dependence in the joint distribution of stock returns. Hu [32] used copulas to study the structure of dependence across financial markets. Miller and Liu [58] proposed methods for recovery of multivariate joint distributions and copulas from limited information using entropy and other information theoretic concepts.

The multivariate measures of dependence and the copula-based approaches to dependence modeling are two interrelated parts of the study of joint distributions of r.v.'s in mathematical statistics and probability theory. A problem of fundamental importance in the field is to determine a relationship between a multivariate cumulative distribution function (cdf) and its lower dimensional margins and to measure degrees of dependence that correspond to particular classes of joint cdf's. The problem is closely related to the problem of characterizing the joint distribution by conditional distributions (see Gouriéroux and Monfort [28]). Remarkable advances have been made in the latter research area in recent years in statistics and probability literature (see, e.g., papers in Dall'Aglio, Kotz and Salinetti [13], Beneš and Štěpán [6] and the monographs by Joe [44], Nelsen [60] and Mari and Kotz [56]).

Motivated by the recent surge in the interest in the study and application of dependence measures and related concepts to account for the complexity in problems in statistics, economics, finance and risk management, this paper provides the first characterizations of joint distributions and copulas for multivariate vectors. These characterizations represent joint distributions of dependent r.v.'s and their copulas as sums of $U$-statistics in independent r.v.'s. We use these characterizations to introduce a unified approach to modeling multivariate dependence and provide new results concerning convergence of multidimensional statistics of time series. The results provide a device for reducing the analysis of convergence of multidimensional statistics of time series to the study of convergence of the measures of intertemporal dependence in the time series (e.g., the multivariate Pearson coefficient, the relative entropy, the multivariate divergence measures, the mean information for discrimination between the dependence and independence, the generalized Tsallis entropy and the Hellinger distance). Furthermore, they allow one to reduce the problems of the study of convergence of statistics of intertemporally dependent time series to the study of convergence of corresponding statistics in the case of intertemporally independent time series. That is, the characterizations for copulas obtained in the paper imply results which associate with each set of arbitrarily dependent r.v.'s a sum of $U$-statistics in independent r.v.'s with canonical kernels. Thus, they allow one to reduce problems for dependent r.v.'s to well-studied objects and to transfer results known for independent r.v.'s and $U$-statistics to the case of arbitrary dependence (see, e.g., Ibragimov and Sharakhmetov [36-40], Ibragimov, Sharakhmetov and Cecen [41], de la Peña, Ibragimov and Sharakhmetov [16, 17] and references therein for general moment inequalities for sums of $U$-statistics and their particular important cases, sums of r.v.'s and multilinear forms, and Ibragimov and Phillips [35] for a new and conceptually simple method for obtaining weak convergence of multilinear forms, $U$-statistics and their non-linear analogues to stochastic integrals based on general asymptotic theory for semimartingales and for applications of the method in a wide range of linear and non-linear time series models).

As a corollary of the results for copulas, we obtain new complete characterizations of important classes of dependent r.v.'s that give, in particular, methods for constructing new copulas and modeling various dependence structures. The results in the paper provide, among others, complete positive answers to the problems raised by Kotz and Seeger [47] concerning characterizations of density weighting functions (d.w.f.) of dependent r.v.'s, existence of methods for constructing d.w.f.'s, and derivation of d.w.f.'s for a given model of dependence (see also [58] for a discussion of d.w.f.'s).

Along the way, a general methodology (of intrinsic interest within and outside probability theory, economics and finance) is developed for analyzing key measures of dependence among r.v.'s. Using the methodology, we obtain sharp decoupling inequalities for comparing the expectations of arbitrary (integrable) functions of dependent variables to their corresponding counterparts with independent variables through the inclusion of multivariate dependence measures.

On the methodological side, the paper shows how the results in theory of $U$-statistics, including inversion formulas for these objects that provide the main tools for the argument for representations in this paper (see the proof of Theorem 1), can be used in the study of joint distributions, copulas and dependence.

The paper is organized as follows. Sections 2 and 3 contain the results on general characterizations of copulas and joint distributions of dependent r.v.'s. Section 4 presents the results on characterizations of dependence based on $U$-statistics in independent r.v.'s. In Sections 5 and 6, we apply the results for copulas and joint distributions to characterize different classes of dependent r.v.'s. Section 7 contains the results on reduction of the analysis of convergence of multidimensional statistics of time series to the study of convergence of the measures of intertemporal dependence in time series as well as the results on sharp decoupling inequalities for dependent r.v.'s. The proofs of the results obtained in the paper are in the Appendix.

## 2. General characterizations of joint distributions of arbitrarily dependent random variables

In the present section, we obtain explicit general representations for joint distributions of arbitrarily dependent r.v.'s absolutely continuous with respect to products of marginal distributions. Let $F_k : \mathbf{R} \to [0,1]$, $k = 1, \ldots, n$, be one-dimensional cdf's and let $\xi_1, \ldots, \xi_n$ be independent r.v.'s on some probability space $(\Omega, \Im, P)$ with $P(\xi_k \leq x_k) = F_k(x_k)$, $x_k \in \mathbf{R}$, $k = 1, \ldots, n$ (we formulate the results for the case of right-continuous cdf's; however, completely similar results hold in the left-continuous case).

In what follows, $F(x_1, \dots, x_n)$, $x_i \in \mathbf{R}$, $i = 1, \dots, n$, stands for a function satisfying the following conditions:

(a) $F(x_1, \dots, x_n) = P(X_1 \le x_1, \dots, X_n \le x_n)$ for some r.v.'s $X_1, \dots, X_n$ on a probability space $(\Omega, \Im, P)$;

(b) the one-dimensional marginal cdf's of $F$ are $F_1, \dots, F_n$;

(c) $F$ is absolutely continuous with respect to $dF(x_1) \cdots dF_n(x_n)$ in the sense that there exists a Borel function $G : \mathbf{R}^n \to [0, \infty)$ such that

$$F(x_1, \dots, x_n) = \int_{-\infty}^{x_1} \cdots \int_{-\infty}^{x_n} G(t_1, \dots, t_n) dF_1(t_1) \cdots dF_n(t_n).$$

As usual, throughout the paper, we denote $G$ in (c) by $\frac{dF}{dF_1 \cdots dF_n}$. In addition, $F(x_{j_1}, \dots, x_{j_k})$, $1 \le j_1 < \cdots < j_k \le n$, $k = 2, \dots, n$, stands for the $k$-dimensional marginal cdf of $F(x_1, \dots, x_n)$. Also, in what follows, if not stated otherwise, $\frac{dF(x_{j_1}, \dots, x_{j_k})}{dF_{j_1} \cdots dF_{j_k}}$, $1 \le j_1 < \cdots < j_k \le n$, $k = 2, \dots, n$, is to be taken to be 1 if at least one $x_{j_1}, \dots, x_{j_k}$ is not a point of increase of the corresponding $F_{j_1}, \dots, F_{j_k}$ (that is, if $(x_{j_1}, \dots, x_{j_k})$ is not in the support of $dF_{j_1} \cdots dF_{j_k}$).

Throughout the paper, the functions $g$ appearing in the representations obtained are assumed to be Borel measurable.

**Theorem 2.1.** *A function $F : \mathbf{R}^n \to [0, 1]$ is a joint cdf with one-dimensional marginal cdf's $F_k(x_k)$, $x_k \in \mathbf{R}$, $k = 1, \dots, n$, absolutely continuous with respect to the product of marginal cdf's $\prod_{k=1}^n F_k(x_k)$, if and only if there exist functions $g_{i_1, \dots, i_c} : \mathbf{R}^c \to \mathbf{R}$, $1 \le i_1 < \cdots < i_c \le n$, $c = 2, \dots, n$, satisfying conditions*

A1 (integrability):

$$E|g_{i_1, \dots, i_c}(\xi_{i_1}, \dots, \xi_{i_c})| < \infty,$$

A2 (degeneracy)*:*

$$E(g_{i_1, \dots, i_c}(\xi_{i_1}, \dots, \xi_{i_{k-1}}, \xi_{i_k}, \xi_{i_{k+1}}, \dots, \xi_{i_c}) | \xi_{i_1}, \dots, \xi_{i_{k-1}}, \xi_{i_{k+1}}, \dots, \xi_{i_c}) =$$

$$\int_{-\infty}^{\infty} g_{i_1, \dots, i_c}(\xi_{i_1}, \dots, \xi_{i_{k-1}}, x_{i_k}, \xi_{i_{k+1}}, \dots, \xi_{i_c}) dF_{i_k}(x_{i_k}) = 0, \ (a.s.)$$

$1 \le i_1 < \cdots < i_c \le n$, $k = 1, 2, \dots, c$, $c = 2, \dots, n$,

A3 (positive definiteness)*:*

$$U_n(\xi_1, \dots, \xi_n) \equiv \sum_{c=2}^{n} \sum_{1 \le i_1 < \cdots < i_c \le n} g_{i_1, \dots, i_c}(\xi_{i_1}, \dots, \xi_{i_c}) \ge -1 \ (a.s.)$$

*and such that the following representation holds for $F$ :*

$$F(x_1, \dots, x_n) = \int_{-\infty}^{x_1} \cdots \int_{-\infty}^{x_n} (1 + U_n(t_1, \dots, t_n)) \prod_{i=1}^{n} dF_i(t_i). \qquad (2.1)$$

*Moreover, $g_{i_1, \dots, i_c}(\xi_{i_1}, \dots, \xi_{i_c}) = f_{i_1, \dots, i_c}(\xi_{i_1}, \dots, \xi_{i_c})$ (a.s.), $1 \le i_1 < \cdots < i_c \le n$, $c = 2, \dots, n$, where*

$$f_{i_1, \dots, i_c}(x_{i_1}, \dots, x_{i_c}) = \sum_{k=2}^{c} (-1)^{c-k} \sum_{1 \le j_1 < \cdots < j_k \in \{i_1, \dots, i_c\}} \left( \frac{dF(x_{j_1}, \dots, x_{j_k})}{dF_{j_1} \cdots dF_{j_k}} - 1 \right).$$

**Remark 2.1.** It is not difficult to see that if r.v.'s $X_1, \ldots, X_n$ have a joint cdf given by (2.1) then the r.v.'s $X_{j_1}, \ldots, X_{j_k}$, $1 \leq j_1 < \cdots < j_k \leq n$, $k = 2, \ldots, n$, have the joint cdf

$$F(x_{j_1}, \ldots, x_{j_k})$$
$$= \int_{-\infty}^{x_{j_1}} \cdots \int_{-\infty}^{x_{j_k}} (1 + \sum_{c=2}^{n} \sum_{\{i_1 < \ldots < i_c\} \in B_k} g_{i_1,\ldots,i_c}(t_{i_1}, \ldots, t_{i_c})) \prod_{i=1}^{k} dF_{j_i}(t_{j_i})$$

with the same functions $g_{i_1,\ldots,i_c}$, and where $B_k = \{j_1, \ldots, j_k\}$.

Theorem 2.1 can be equivalently formulated as in the following remark.

**Remark 2.2.** A function $F : \mathbf{R}^n \to [0,1]$ is a joint cdf with the one-dimensional marginal cdf's $F_k(x_k)$, $x_k \in \mathbf{R}$, $k = 1, \ldots, n$, absolutely continuous with respect to the product of marginal cdf's $\prod_{k=1}^{n} F_k(x_k)$, if and only if there exist functions $g_{i_1,\ldots,i_c} : \mathbf{R}^c \to \mathbf{R}$, $1 \leq i_1 < \cdots < i_c \leq n$, $c = 2, \ldots, n$, satisfying conditions A1–A3 and such that the element of frequency $dF(x_1, \ldots, x_n)$ can be expressed in the form

$$dF(x_1, \ldots, x_n) = \prod_{i=1}^{n} dF_i(x_i)(1 + U_n(x_1, \ldots, x_n)). \tag{2.2}$$

**Remark 2.3.** Sharakhmetov [69] provided proof of (2.2) in the case of density functions (of r.v.'s absolutely continuous with respect to *Lebesgue* measure), with a mention that a similar representation holds for distributions of discrete r.v.'s. The setup considered in this paper includes (among others) the class of vectors of dependent absolutely continuous and discrete r.v.'s as well as vectors of mixtures of absolutely continuous and discrete r.v.'s. Furthermore, our proof easily extends to the case of general Banach spaces, in particular, the spaces $\mathbf{R}^k$.

## 3. Applications to copulas

Let us start with the definition of copulas and the formulation of Sklar's theorem mentioned in the introduction (see e.g., [22] and [60]).

**Definition 3.1.** A function $C : [0,1]^n \to [0,1]$ is called a $n$-dimensional *copula* if it satisfies the following conditions:

1. $C(u_1, \ldots, u_n)$ is increasing in each component $u_i$.
2. $C(u_1, \ldots, u_{k-1}, 0, u_{k+1}, \ldots, u_n) = 0$ for all $u_i \in [0,1]$, $i \neq k$, $k = 1, \ldots, n$.
3. $C(1, \ldots, 1, u_i, 1, \ldots, 1) = u_i$ for all $u_i \in [0,1]$, $i = 1, \ldots, n$.
4. For all $(a_1, \ldots, a_n)$, $(b_1, \ldots, b_n) \in [0,1]^n$ with $a_i \leq b_i$,

$$\sum_{i_1=1}^{2} \cdots \sum_{i_n=1}^{2} (-1)^{i_1+\cdots+i_n} C(x_{1i_1}, \ldots, x_{ni_n}) \geq 0,$$

where $x_{j1} = a_j$ and $x_{j2} = b_j$ for all $j \in \{1, \ldots, n\}$. Equivalently, $C$ is a $n$-dimensional copula if it is a joint cdf of $n$ r.v.'s each of which is uniformly distributed on $[0,1]$.

**Definition 3.2.** A copula $C : [0,1]^n \to [0,1]$ is called *absolutely continuous* if, when considered as a joint cdf, it has a joint density given by $\partial C^n(u_1 \ldots, u_n)/\partial u_1 \cdots \partial u_n$.

**Theorem 3.1** (Sklar [72]). *If $X_1, \dots, X_n$ are random variables defined on a common probability space, with the one-dimensional cdf's $F_{X_k}(x_k) = P(X_k \le x_k)$ and the joint cdf $F_{X_1,\dots,X_n}(x_1, \dots, x_n) = P(X_1 \le x_1, \dots, X_n \le x_n)$, then there exists an n-dimensional copula $C_{X_1,\dots,X_n}(u_1, \dots, u_n)$ such that $F_{X_1,\dots,X_n}(x_1, \dots, x_n) = C_{X_1,\dots,X_n}(F_{X_1}(x_1), \dots, F_{X_n}(x_n))$ for all $x_k \in \mathbf{R}$, $k = 1, \dots, n$.*

The following theorems give analogues of the representations in the previous section for copulas. Let $V_1, \dots, V_n$ denote independent r.v.'s uniformly distributed on $[0, 1]$.

**Theorem 3.2.** *A function $C : [0, 1]^n \to [0, 1]$ is an absolutely continuous n-dimensional copula if and only if there exist functions $\tilde{g}_{i_1,\dots,i_c} : \mathbf{R}^c \to \mathbf{R}$, $1 \le i_1 < \cdots < i_c \le n$, $c = 2, ..., n$, satisfying the conditions*
A4 (integrability)*:*

$$\int_0^1 \cdots \int_0^1 |\tilde{g}_{i_1,\dots,i_c}(t_{i_1}, \dots, t_{i_c})| dt_{i_1} \cdots dt_{i_c} < \infty,$$

A5 (degeneracy)*:*

$$E(\tilde{g}_{i_1,\dots,i_c}(V_{i_1}, \dots, V_{i_{k-1}}, V_{i_k}, V_{i_{k+1}}, \dots, V_{i_c}) | V_{i_1}, \dots, V_{i_{k-1}}, V_{i_{k+1}}, \dots, V_{i_c}) \\ = \int_0^1 \tilde{g}_{i_1,\dots,i_c}(V_{i_1}, \dots, V_{i_{k-1}}, t_{i_k}, V_{i_{k+1}}, \dots, V_{i_c}) dt_{i_k} = 0 \,(a.s.),$$

$1 \le i_1 < \cdots < i_c \le n$, $k = 1, 2, \dots, c$, $c = 2, \dots, n$,

A6 (positive definiteness):

$$\tilde{U}_n(V_1, \dots, V_n) \equiv \sum_{c=2}^{n} \sum_{1 \le i_1 < \cdots < i_c \le n} \tilde{g}_{i_1,\dots,i_c}(V_{i_1}, \dots, V_{i_c}) \ge -1 \,(a.s.)$$

*and such that*

$$C(u_1, \dots, u_n) = \int_0^{u_1} \cdots \int_0^{u_n} (1 + \tilde{U}_n(t_1, \dots, t_n)) \prod_{i=1}^{n} dt_i. \tag{3.1}$$

Theorem 3.2 and Sklar's theorem formulated above imply the following representation for a joint distribution of r.v.'s.

**Theorem 3.3.** *A function $F : \mathbf{R}^n \to [0, 1]$ is a joint cdf with one-dimensional marginal cdf's $F_k(x_k)$, $x_k \in \mathbf{R}$, $k = 1, \dots, n$, absolutely continuous with respect to the product of marginal cdf's $\prod_{k=1}^n F_k(x_k)$ if and only if there exist functions $\tilde{g}_{i_1,\dots,i_c} : [0, 1]^c \to \mathbf{R}$, $1 \le i_1 < \cdots < i_c \le n$, $c = 2, \dots, n$, satisfying conditions* A4–A6 *and such that the following representation holds for $F$:*

$$F(x_1, \dots, x_n) = \int_0^{F_1(x_1)} \cdots \int_0^{F_n(x_n)} (1 + \tilde{U}_n(t_1, \dots, t_n)) \prod_{i=1}^{n} dt_i, \tag{3.2}$$

*or, equivalently, if and only if the element of joint frequency $dF$ can be expressed in the form*

$$dF(x_1, \dots, x_n) = \prod_{i=1}^{n} dF_i(x_i)(1 + \tilde{U}_n(F_1(x_1), \dots, F_n(x_n))).$$

**Remark 3.1.** The functions $g$ and $\tilde{g}$ in Theorems 2.1-3.3 are related in the following way: $g_{i_1,\dots,i_c}(x_{i_1},\dots,x_{i_c}) = \tilde{g}_{i_1,\dots,i_c}(F_{i_1}(x_{i_1}),\dots,F_{i_c}(x_{i_c}))$.

Theorems 2.1–3.3 provide a general device for constructing multivariate copulas and distributions. E.g., taking in (3.1) and (3.2) $n = 2$, $\tilde{g}_{1,2}(t_1,t_2) = \alpha(1-2t_1)(1-2t_2)$, $\alpha \in [-1,1]$, we get the family of bivariate Eyraud–Farlie–Gumbel–Morgenstern copulas $C_\alpha(u_1,u_2) = u_1u_2(1+\alpha(1-u_1)(1-u_2))$ and corresponding distributions $F_\alpha(x_1,x_2) = F_1(x_1)F_2(x_2)(1+\alpha(1-F_1(x_1))(1-F_2(x_2))$. More generally, taking $\tilde{g}_{i_1,\dots,i_c}(t_{i_1},\dots,t_{i_c}) = 0$, $1 \le i_1 < \cdots < i_c \le n$, $c = 2,\dots,n-1$, $\tilde{g}_{1,2,\dots,n}(t_1,t_2,\dots,t_n) = \alpha(1-2t_1)(1-2t_2)\cdots(1-2t_n)$, we obtain the multivariate Eyraud–Farlie–Gumbel–Morgenstern copulas $C_\alpha(u_1,u_2,\dots,u_n) = \prod_{i=1}^n u_i(1+\alpha\prod_{i=1}^n(1-u_i))$ and corresponding multivariate cdf's $F_\alpha(x_1,x_2,\dots,x_n) = \prod_{i=1}^n F_i(x_i)(1+\alpha\prod_{i=1}^n(1-F_i(x_i)))$.

Let $\alpha_{i_1,\dots,i_c} \in \mathbf{R}$ be constants such that $\sum_{c=2}^n \sum_{1\le i_1<\cdots<i_c\le n} \alpha_{i_1,\dots,i_c} \times \delta_{i_1}\cdots\delta_{i_c} \ge -1$ for all $\delta_i \in \{0,1\}$, $i = 1,\dots,n$. The choice $\tilde{g}_{i_1,\dots,i_c}(t_{i_1},\dots,t_{i_c}) = \alpha_{i_1,\dots,i_c}(1-2t_{i_1})(1-2t_{i_2})\cdots(1-2t_{i_c})$, $1 \le i_1 < \cdots < i_c \le n$, $c = 2,\dots,n$, gives the following generalized multivariate Eyraud–Farlie–Gumbel–Morgenstern copulas (see Johnson and Kotz [45] and Cambanis [10]):

$$C(u_1,\dots,u_n) = \prod_{k=1}^n u_k(1+\sum_{c=2}^n \sum_{1\le i_1<\cdots<i_c\le n} \alpha_{i_1,\dots,i_c}(1-u_{i_k})) \tag{3.3}$$

and the corresponding cdf's

$$F(x_1,\dots,x_n) = \prod_{i=1}^n F_i(x_i)(1+\sum_{c=2}^n \sum_{1\le i_1<\cdots<i_c\le n} \alpha_{i_1,\dots,i_c}(1-F_{i_k}(x_{i_k}))).$$

The importance of the generalized Eyraud–Farlie–Gumbel–Morgenstern copulas and cdf's stems, in particular, from the fact that, as shown in Sharakhmetov and Ibragimov [70], they completely characterize joint distributions of two-valued r.v.'s.

Taking $n = 2$, $\tilde{g}_{1,2}(t_1,t_2) = \theta c(t_1,t_2)$, where $c$ is a continuous function on the unit square $[0,1]^2$ satisfying the properties $\int_0^1 c(t_1,t_2)dt_1 = \int_0^1 c(t_1,t_2)dt_2 = 0$, $1+\theta c(t_1,t_2) \ge 0$ for all $0 \le t_1,t_2 \le 1$, one obtains the class of bivariate densities studied by Rüschendorf [67] and Long and Krzysztofowicz [53] (see also [56, pp. 73–78]) $f(x_1,x_2) = f_1(x_1)f_2(x_2)(1+\theta c(F_1(x_1),F_2(x_2)))$ with the covariance characteristic $c$ and the covariance scalar $\theta$. Furthermore, from Theorems 2.1–3.3 it follows that this representation in fact holds for an arbitrary density function and the function $\theta c(t_1,t_2)$ is unique.

## 4. From dependence to independence through $U$-statistics

Denote by $\mathcal{G}_n$ the class of sums of $U$-statistics of the form

$$U_n(\xi_1,\dots,\xi_n) = \sum_{c=2}^n \sum_{1\le i_1<\cdots<i_c\le n} g_{i_1,\dots,i_c}(\xi_{i_1},\dots,\xi_{i_c}), \tag{4.1}$$

where the functions $g_{i_1,\dots,i_c}$, $1 \le i_1 < \cdots < i_c \le n$, $c = 2,\dots,n$, satisfy conditions A1–A3, and, as before, $\xi_1,\dots,\xi_n$ are independent r.v.'s with cdf's $F_k(x_k)$, $x_k \in \mathbf{R}$, $k = 1,\dots,n$.

The following theorem puts into correspondence to any set of arbitrarily dependent r.v.'s a sum of $U$-statistics in independent r.v.'s with canonical kernels. This allows one to reduce problems for dependent r.v.'s to well-studied objects and to transfer results known for independent r.v.'s and $U$-statistics to the case of arbitrary dependence. In what follows, the joint distributions considered are assumed to be absolutely continuous with respect to the product of the marginal distributions $\prod_{k=1}^n F_k(x_k)$.

**Theorem 4.1.** *The r.v.'s $X_1, \ldots, X_n$ have one-dimensional cdf's $F_k(x_k)$, $x_k \in \mathbf{R}$, $k = 1, \ldots, n$, if and only if there exists $U_n \in \mathcal{G}_n$ such that for any Borel measurable function $f : \mathbf{R}^n \to \mathbf{R}$ for which the expectations exist*

$$Ef(X_1, \ldots, X_n) = Ef(\xi_1, \ldots, \xi_n)(1 + U_n(\xi_1, \ldots, \xi_n)). \tag{4.2}$$

Note that the above Theorem 4.1 holds for complex-valued functions $f$ as well as for real-valued ones. That is, letting $f(x_1, \ldots, x_n) = \exp(i \sum_{k=1}^n t_k x_k)$, $t_k \in \mathbf{R}$, $k = 1, \ldots, n$, one gets the following representation for the joint characteristic function of the r.v.'s $X_1, \ldots, X_n$ :

$$E \exp\left(i \sum_{k=1}^n t_k X_k\right) = E \exp\left(i \sum_{k=1}^n t_k \xi_k\right) + E \exp\left(i \sum_{k=1}^n t_k \xi_k\right) U_n(\xi_1, \ldots, \xi_n).$$

## 5. Characterizations of classes of dependent random variables

The following Theorems 5.1–5.8 give characterizations of different classes of dependent r.v.'s in terms of functions $g$ that appear in the representations for joint distributions obtained in Section 2. Completely similar results hold for the functions $\tilde{g}$ that enter corresponding representations for copulas in Section 3.

**Theorem 5.1.** *The r.v.'s $X_1, \ldots, X_n$ with one-dimensional cdf's $F_k(x_k)$, $x_k \in \mathbf{R}$, $k = 1, \ldots, n$, are independent if and only if the functions $g_{i_1,\ldots,i_c}$ in representations* (2.1) *and* (2.2) *satisfy the conditions*

$$g_{i_1,\ldots,i_c}(\xi_{i_1}, \ldots, \xi_{i_c}) = 0 \text{ (a.s.)},\ 1 \le i_1 < \cdots < i_c \le n, c = 2, \ldots, n.$$

**Theorem 5.2.** *A sequence of r.v.'s $\{X_n\}$ is strictly stationary if and only if the functions $g_{i_1,\ldots,i_c}$ in representations* (2.1) *and* (2.2) *for any finite-dimensional distribution (see Remark* 2.1*) satisfy the conditions*

$$g_{i_1+h,\ldots,i_c+h}(\xi_{i_1}, \ldots, \xi_{i_c}) = g_{i_1,\ldots,i_c}(\xi_{i_1}, \ldots, \xi_{i_c}) \text{ (a.s.) } 1 \le i_1 < \cdots < i_c \le n,$$

$$c = 2, 3, \ldots, h = 0, 1, \ldots$$

**Theorem 5.3.** *A sequence of r.v.'s $\{X_n\}$ with $EX_k = 0$, $EX_k^2 < \infty$, $k = 1, 2, \ldots$, is weakly stationary if and only if the functions $g$ in representations* (2.1) *and* (2.2) *for any finite-dimensional distribution have the property that the function $h(s, t) = E\xi_s \xi_t g_{st}(\xi_s, \xi_t)$, depends only on $|t - s|$, $t, s = 1, 2, \ldots$*

**Definition 5.1.** The r.v.'s $X_1, \ldots, X_n$ with $EX_i = 0$, $i = 1, \ldots, n$, are called *orthogonal* if $EX_i X_j = 0$ for all $1 \le i < j \le n$.

**Theorem 5.4.** *The r.v.'s $X_1, \ldots, X_n$ with $EX_k = 0$, $k = 1, \ldots, n$, are orthogonal if and only if the functions $g$ in representations* (2.1) *and* (2.2) *satisfy the conditions $E\xi_i \xi_j g_{ij}(\xi_i, \xi_j) = 0$, $1 \le i < j \le n$.*

**Definition 5.2.** The r.v.'s $X_1, \ldots, X_n$ are called exchangeable if all $n!$ permutations $(X_{\pi(1)}, \ldots, X_{\pi(n)})$ of the r.v.'s have the same joint distributions.

**Theorem 5.5.** *The identically distributed r.v.'s $X_1, \ldots, X_n$ are exchangeable if and only if the functions $g_{i_1,\ldots,i_c}$ in representations* (2.1) *and* (2.2) *satisfy the conditions $g_{i_1,\ldots,i_c}(\xi_{i_1}, \ldots, \xi_{i_c}) = g_{i_{\pi(1)},\ldots,i_{\pi(c)}}(\xi_{i_{\pi(1)}}, \ldots, \xi_{i_{\pi(i_c)}})$ (a.s.) for all $1 \le i_1 < \cdots < i_c \le n$, $c = 2, \ldots, n$, and all permutations $\pi$ of the set $\{1, \ldots, n\}$.*

**Definition 5.3.** The r.v.'s $X_1, \ldots, X_n$ are called *$m$-dependent* $(1 \le m \le n)$ if any two vectors $(X_{j_1}, X_{j_2}, \ldots, X_{j_{a-1}}, X_{j_a})$ and $(X_{j_{a+1}}, X_{j_{a+2}}, \ldots, X_{j_{l-1}}, X_{j_l})$, where $1 \le j_1 < \cdots < j_a < \cdots < j_l \le n$, $a = 1, 2, \ldots, l-1$, $l = 2, \ldots, n$, $j_{a+1} - j_a \ge m$, are independent.

**Theorem 5.6.** *The r.v.'s $X_1, \ldots, X_n$ are $m$-dependent if and only if the functions $g$ in representations* (2.1) *and* (2.2) *satisfy the conditions $g_{i_1,\ldots,i_k,i_{k+1},\ldots,i_c}(\xi_{i_1}, \ldots, \xi_{i_k}, \xi_{i_{k+1}}, \ldots, \xi_{i_c}) = g_{i_1,\ldots,i_k}(\xi_{i_1}, \ldots, \xi_{i_k}) g_{i_{k+1},\ldots,i_c}(\xi_{i_{k+1}}, \ldots, \xi_{i_c})$ for all $1 \le i_1 < \cdots < i_k < i_{k+1} < \cdots < i_c \le n$, $i_{k+1} - i_k \ge m$, $k = 1, \ldots, c-1$, $c = 2, \ldots, n$.*

**Definition 5.4.** The r.v.'s $X_1, \ldots, X_n$ form a *multiplicative system* of order $\alpha \in \mathbf{N}$ (shortly, $MS(\alpha)$) if $E|X_j|^\alpha < \infty$, $j = 1, \ldots, n$, and $E\prod_{j=1}^n X_j^{\alpha_j} = \prod_{j=1}^n EX_j^{\alpha_j}$ for any $\alpha_j \in \{0, 1, \ldots, \alpha\}$, $j = 1, \ldots, n$.

The systems $MS(1)$ and $MS(2)$ under the names multiplicative and strongly multiplicative systems, respectively, were introduced by Alexits [2]. Multiplicative systems of an arbitrary order were considered, e.g., by Kwapień [48] and Sharakhmetov [68]. Examples of the multiplicative systems $MS(1)$ are given, besides independent r.v.'s, by the lacunary trigonometric systems $\{\cos 2\pi n_k x, \sin 2\pi n_k x, k = 1, 2, \ldots\}$ on the interval $[0, 1]$ with the Lebesgue measure for $n_{k+1}/n_k \ge 2$ important in Fourier analysis of time series and also by such important classes of dependent r.v.'s as martingale-difference sequences. Examples of strongly multiplicative systems (that is, the systems $MS(2)$) are given by the lacunary trigonometric systems for $n_{k+1}/n_k \ge 3$ and martingale-difference sequences $X_1, \ldots, X_n$ satisfying the conditions $E(X_n^2 | X_1, \ldots, X_{n-1}) = b_n^2 \in \mathbf{R}$, $n = 1, 2, \ldots$. Examples of the systems $MS(\alpha)$ include, for instance, the lacunary trigonometric systems with large lacunas, that is, with $n_{k+1}/n_k \ge \alpha + 1$ and also $\epsilon$-independent and asymptotically independent r.v.'s introduced by Zolotarev [78].

**Theorem 5.7.** *The r.v.'s $X_1, \ldots, X_n$ form a multiplicative system of order $\alpha$ if and only if the functions $g_{i_1,\ldots,i_c}$ in representations* (2.1) *and* (2.2) *satisfy the conditions $E\xi_{i_1}^{\alpha_{i_1}} \cdots \xi_{i_c}^{\alpha_{i_c}} g_{i_1,\ldots,i_c}(\xi_{i_1}, \ldots, \xi_{i_c}) = 0$, $1 \le i_1 < \cdots < i_c \le n$, $c = 2, \ldots, n$, $\alpha_j \in \{0, 1, \ldots, \alpha\}$, $j = 1, \ldots, n$.*

**Definition 5.5.** The r.v.'s $X_1, \ldots, X_n$ are called *$r$-independent* $(2 \le r < n)$ if any $r$ of them of are jointly independent.

**Theorem 5.8.** *The r.v.'s $X_1, \ldots, X_n$ are $r$-independent if and only if the functions $g_{i_1,\ldots,i_c}$ in representations* (2.1) *and* (2.2) *satisfy the conditions $g_{i_1,\ldots,i_c}(\xi_{i_1}, \ldots, \xi_{i_c}) = 0$ (a.s.), $1 \le i_1 < \cdots < i_c \le n$, $c = 2, \ldots, r$.*

**Remark 5.1.** Let $F_1(x), \ldots, F_n(x)$ be arbitrary one-dimensional distribution functions, $\alpha_1, ..., \alpha_n \in (-1, 1) \setminus \{0\}$, $\sum_{i=1}^n |\alpha_i| \le 1$. Taking $g_{i_1,\ldots,i_c}(t_{i_1}, \ldots, t_{i_c}) = 0$, $1 \le i_1 < \cdots < i_c \le n$, $c = 2, \ldots, n$, $c \ne r+1$, $g_{i_1,\ldots,i_{r+1}}(t_{i_1}, \ldots, t_{i_{r+1}}) = \frac{\alpha_1 \cdots \alpha_n}{\alpha_{i_1} \cdots \alpha_{i_{r+1}}}((k+1)t_{i_1}^k - (k+2)t_{i_1}^{k+1}) \cdots ((k+1)t_{i_c}^k - (k+2)t_{i_c}^{k+1})$, $k = 0, 1, 2, \ldots$, in Theorem 3.3, we obtain the following extensions of the examples of $r$-independent

r.v.'s obtained by Wang [77]: For $k = 0, 1, 2, \ldots$,

$$F(x_1, \ldots, x_n) = \prod_{i=1}^{n} F_i(x_i) \Bigg( 1 + \sum_{1 \le i_1 < \ldots < i_{r+1} \le n} \frac{\alpha_1 \ldots \alpha_n}{\alpha_{i_1} \ldots \alpha_{i_r+1}} (F_{i_1}^k(x_{i_1}) - F_{i_1}^{k+1}(x_{i_1})) \cdots \times (F_{i_{r+1}}^k(x_{i_{r+1}}) - F_{i_{r+1}}^{k+1}(x_{i_{r+1}})) \Bigg),$$

(Wang's examples are with $k = 0$).

## 6. Further applications: a structural property of multiplicative systems

The following theorem shows that r.v.'s forming a multiplicative system of order $\alpha$ and taking not more than $\alpha+1$ values are jointly independent. Let $card(A_i)$ denote the number of elements in (finite) sets $A_i$.

**Theorem 6.1.** *Let $\alpha \in \mathbf{N}$, and let $A_i$, $i = 1, \ldots, n$, be sets of real numbers such that $card(A_i) \le \alpha + 1$, $i = 1, \ldots, n$. The r.v.'s $X_1, \ldots, X_n$ taking values in $A_1, \ldots, A_n$, respectively, form a multiplicative system of order $\alpha$ if and only if they are jointly independent.*

**Remark 6.1.** From Theorem 6.1 with $\alpha = 1$ the following result obtained in [70] follows: A sequence of r.v.'s $\{X_n\}$ on a probability space $(\Omega, \Im, P)$ assuming two values is a martingale-difference with respect to an increasing sequence of $\sigma$-algebras $\Im_0 = (\Omega, \emptyset) \subseteq \Im_1 \subseteq \cdots \subseteq \Im$ if and only if the r.v.'s $\{X_n\}$ are jointly independent. In addition, we obtain that if a sequence of r.v.'s $\{X_n\}$ assuming three values is a martingale-difference with respect to $(\Im_n)$ such that $E(X_n^2 | \Im_{n-1}) = b_n^2 \in \mathbf{R}$, then the r.v.'s are jointly independent.

## 7. Measures of dependence and sharp moment and probability inequalities for dependent random variables

In this section, we apply the results from Section 2 to study properties of different measures of dependence and convergence of multidimensional statistics of time series. We obtain results that allow one to reduce the analysis of convergence of statistics of time series to the study of convergence of the measures of intertemporal dependence in the time series and limit behavior of the statistics in the case of independence. We also prove new sharp complete decoupling inequalities for dependent r.v.'s in terms of their dependence characteristics. The theory of complete decoupling inequalities has experienced an impetus in recent years. The interested reader should consult de la Peña [14], de la Peña and Giné [15] and de la Peña and Lai [18] (a survey) for more on the subject. Let $X_1, \ldots, X_n$ be r.v.'s with one-dimensional cdf's $F_k(x_k)$, $k = 1, \ldots, n$, and joint cdf $F(x_1, \ldots, x_n)$. Recall that $G(x_1, \ldots, x_n) = dF(x_1, \ldots, x_n) / \prod_{i=1}^n dF_i$ and consider the following measures of dependence for the r.v.'s $X_1, \ldots, X_n$ :

$$\begin{aligned} \phi^2_{X_1,\ldots,X_n} &= \int_{-\infty}^{\infty} \cdots \int_{-\infty}^{\infty} G(x_1, \ldots, x_n) dF(x_1, \ldots, x_n) - 1 \\ &= \int_{-\infty}^{\infty} \cdots \int_{-\infty}^{\infty} G^2(x_1, \ldots, x_n) \prod_{i=1}^n dF_i(x_i) - 1 \end{aligned}$$

(multivariate analog of Pearson's $\phi^2$ coefficient), and

$$\delta_{X_1,\dots,X_n} = \int_{-\infty}^{\infty} \cdots \int_{-\infty}^{\infty} \log(G(x_1,\dots,x_n))dF(x_1,...,x_n)$$

(relative entropy), where the integral signs are in the sense of Lebesgue–Stieltjes and $G(x_1,\dots,x_n)$ is taken to be 1 if $(x_1,\dots,x_n)$ is not in the support of $dF_1\cdots dF_n$. In the case of absolutely continuous r.v.'s $X_1,\dots,X_n$ the measures $\delta_{X_1,\dots,X_n}$ and $\phi^2_{X_1,\dots,X_n}$ were introduced by Joe [42, 43]. In the case of two r.v.'s $X_1$ and $X_2$ the measure $\phi^2_{X_1,X_2}$ was introduced by Pearson [63] and was studied, among others, by Lancaster [49–51]. In the bivariate case, the measure $\delta_{X_1,X_2}$ is commonly known as Shannon or Kullback–Leibler mutual information between $X_1$ and $X_2$. It should be noted (see [43]) that if $(X_1,\dots,X_n)' \sim N(\mu,\Sigma)$, then $\phi^2_{X_1,\dots,X_n} = |R(2I_n - R)|^{-1/2} - 1$, where $I_n$ is the $n\times n$ identity matrix, provided that the correlation matrix $R$ corresponding to $\Sigma$ has the maximum eigenvalue of less than 2 and is infinite otherwise ($|A|$ denotes the determinant of a matrix $A$). In addition to that, if in the above case $\mathrm{diag}(\Sigma) = (\sigma_1^2,\dots,\sigma_n^2)$, then $\delta_{X_1,\dots,X_n} = -.5\log(|\Sigma|/\prod_{i=1}^n \sigma_i^2)$. In the case of two normal r.v.'s $X_1$ and $X_2$ with the correlation coefficient $\rho$, $(\phi^2_{X_1,X_2}/(1+\phi^2_{X_1,X_2}))^{1/2} = (1-\exp(-2\delta_{X_1,X_2}))^{1/2} = |\rho|$.

The multivariate Pearson's $\phi^2$ coefficient and the relative entropy are particular cases of multivariate divergence measures $D^{\psi}_{X_1,\dots,X_n} = \int_{-\infty}^{\infty}\cdots\int_{-\infty}^{\infty}\psi(G(x_1,\dots,x_n))\prod_{i=1}^n dF_i(x_i)$, where $\psi$ is a strictly convex function on $\mathbf{R}$ satisfying $\psi(1)=0$ and $G(x_1,\dots,x_n)$ is to be taken to be 1 if at least one $x_1,\dots,x_n$ is not a point of increase of the corresponding $F_1,\dots,F_n$. Bivariate divergence measures were considered, e.g., by Ali and Silvey [3] and Joe [43]. The multivariate Pearson's $\phi^2$ corresponds to $\psi(x) = x^2 - 1$ and the relative entropy is obtained with $\psi(x) = x\log x$.

A class of measures of dependence closely related to the multivariate divergence measures is the class of generalized entropies introduced by Tsallis [75] in the study of multifractals and generalizations of Boltzmann–Gibbs statistics (see also [24, 26, 27])

$$\rho^{(q)}_{X_1,\dots,X_n} = \frac{1}{1-q}(1 - \int_{-\infty}^{\infty}\cdots\int_{-\infty}^{\infty} G^{1-q}(x_1,\dots,x_n))\prod_{i=1}^n dF_i(x_i),$$

where $q$ is the entropic index. In the limiting case $q\to 1$, the discrepancy measure $\rho^{(q)}$ becomes the relative entropy $\delta_{X_1,\dots,X_n}$ and in the case $q \to 1/2$ it becomes the scaled squared Hellinger distance between $dF$ and $dF_1\cdots dF_n$

$$\rho^{(1/2)}_{X_1,\dots,X_n} = \frac{1}{2}(1 - \int_{-\infty}^{\infty}\cdots\int_{-\infty}^{\infty} G^{1/2}(x_1,\dots,x_n))\prod_{i=1}^n dF_i(x_i)) = 2\mathcal{H}^2_{X_1,\dots,X_n}$$

($\mathcal{H}_{X_1,\dots,X_n}$ stands for the Hellinger distance). The generalized entropy has the form of the multivariate divergence measures $D^{\psi}_{X_1,\dots,X_n}$ with $\psi(x) = (1/(1-q))(1-x^{1-q})$.

In the terminology of information theory (see, e.g., Akaike [1]) the multivariate analog of Pearson coefficient, the relative entropy and, more generally, the multivariate divergence measures represent the mean amount of information for discrimination between the density $f$ of dependent sample and the density of the sample of independent r.v.'s with the same marginals $f_0 = \prod_{k=1}^n f_k(x_k)$ when the actual distribution is dependent $I(f_0, f; \Phi) = \int \Phi(f(x)/f_0(x))f(x)dx$, where $\Phi$ is a properly chosen function. The multivariate analog of Pearson coefficient is characterized by the relation (below, $f_0$ denotes the density of independent sample and $f$ denotes the density of a dependent sample) $\phi^2 = I(f_0, f; \Phi_1)$, where $\Phi_1(x) = x$; the relative entropy satisfies $\delta = I(f_0, f; \Phi_2)$, where $\Phi_2(x) = \log(x)$; and the multivariate divergence measures satisfy $D^{\psi}_{X_1,\ldots,X_n} = I(f_0, f, \Phi_3)$, where $\Phi_3(x) = \psi(x)/x$.

If $g_{i_1,\ldots,i_c}(x_{i_1}, \ldots, x_{i_c})$ are functions corresponding to Theorem 2.1 and Remark 2.2, then from Theorem 4.1 it follows that the measures $\delta_{X_1,\ldots,X_n}$, $\phi^2_{X_1,\ldots,X_n}$, $D^{\psi}_{X_1,\ldots,X_n}$, $\rho^{(q)}_{X_1,\ldots,X_n}$ (in particular, $2\mathcal{H}^2_{X_1,\ldots,X_n}$ for $q = 1/2$) and $I(f_0, f; \Phi)$ can be written as

$$\begin{aligned} \delta_{X_1,\ldots,X_n} &= E \log\left(1 + U_n(X_1, \ldots, X_n)\right) \\ &= E\left(1 + U_n(\xi_1, \ldots, \xi_n)\right) \log(1 + U_n(\xi_1, \ldots, \xi_n))\,, \end{aligned} \tag{7.1}$$

$$\begin{aligned} \phi^2_{X_1,\ldots,X_n} &= E\left(1+U_n(\xi_1, \ldots, \xi_n)\right)^2 - 1 \\ &= EU_n^2(\xi_1, \ldots, \xi_n) = EU_n(X_1, \ldots, X_n), \end{aligned} \tag{7.2}$$

$$D^{\psi}_{X_1,\ldots,X_n} = E\psi\left(1 + U_n(\xi_1, \ldots, \xi_n)\right), \tag{7.3}$$

$$\rho^{(q)}_{X_1,\ldots,X_n} = (1/(1-q))(1 - E(1 + U_n(\xi_1, \ldots, \xi_n))^q), \tag{7.4}$$

$$2\mathcal{H}^2_{X_1,\ldots,X_n} = 1/2(1 - E(1 + U_n(\xi_1, \ldots, \xi_n))^{1/2}), \tag{7.5}$$

$$I(f_0, f; \Phi) = E\Phi\left(1 + U_n(\xi_1, \ldots, \xi_n)\right)\left(1 + U_n(\xi_1, \ldots, \xi_n)\right), \tag{7.6}$$

where $U_n(x_1, \ldots, x_n)$ is as defined by (4.1).

From (7.2) it follows that the following formula that gives an expansion for $\phi^2_{X_1,\ldots,X_n}$ in terms of the "canonical" functions $g$ holds: $\phi^2_{X_1,\ldots,X_n} = \sum_{c=2}^{n} \sum_{1 \le i_1 < \cdots < i_c \le n} Eg^2_{i_1,\ldots,i_c}(\xi_{i_1}, \ldots, \xi_{i_c})$. In particular, in the case of the r.v.'s $X_1, \ldots, X_n$ with the generalized multivariate Eyraud–Farlie–Gumbel–Morgenstern copulas (3.3) the measure of dependence $\phi^2_{X_1,\ldots,X_n}$ is given by $\phi^2_{X_1,\ldots,X_n} = \sum_{c=2}^{n} \sum_{1 \le i_1 < \cdots < i_c \le n} \alpha^2_{i_1,\ldots,i_c}$.

It is well known that the mutual information between two r.v.'s $X_1$ and $X_2$ is nonnegative (see [12, p. 27]). The multivariate analog of this property for $\delta_{X_1,\ldots,X_n}$ follows from the results obtained by Joe [43]. It is interesting to note that nonnegativity of $\delta_{X_1,\ldots,X_n}$ can be easily obtained from (7.1): since the function $(1 + x)\ln(1 + x)$ is convex in $x \ge 0$, by Jensen inequality we get $\delta_{X_1,\ldots,X_n} \ge (1 + EU_n(\xi_1, \ldots, \xi_n)) \log(1 + EU_n(\xi_1, \ldots, \xi_n)) = 0$. The following theorem gives an inequality between the measures $\delta_{X_1,\ldots,X_n}$ and $\phi^2_{X_1,\ldots,X_n}$ that generalizes and improves the results obtained, in the bivariate case, by Dragomir [20] (see also Mond and Pečarič [59]).

**Theorem 7.1.** *The following inequalities hold*

$$\delta_{X_1,\ldots,X_n} \le \log(1 + \phi^2_{X_1,\ldots,X_n}) \le \phi^2_{X_1,\ldots,X_n}.$$

The results obtained in the previous sections, in particular, Theorem 4.1, provide a device for reducing the analysis of convergence in distribution of double arrays of dependent variables to the study of the convergence in distribution of a decoupled counterpart plus a measure of dependence in the time series. We apply this idea in reduction of the analysis of convergence of multidimensional statistics of time series to the study of convergence of the measures of intertemporal dependence of the time series, including the above multivariate Pearson coefficient $\phi$, the relative entropy $\delta$, the divergence measures $D^{\psi}$ and the mean information for discrimination between the dependence and independence $I(f_0, f; \Phi)$. We obtain the following Theorem 7.2 which deals with the convergence in distribution of $m$-dimensional statistics of time series.

Let $h : \mathbf{R}^m \to \mathbf{R}$ be an arbitrary function of $m$ arguments, $Y$ be some r.v. and let $\psi$ be a convex function increasing on $[1, \infty)$ and decreasing on $(-\infty, 1)$ with $\psi(1) = 0$. In what follows, $\stackrel{\mathcal{D}}{\to}$ represents convergence in distribution. In addition, $\{\xi_i^n\}$ and $\{\xi_t\}$ stand for dependent copies of $\{X_i^n\}$ and $\{X_t\}$.

**Theorem 7.2.** *For the double array* $\{X_i^n\}, i = 1, \dots, n, n = 0, 1, \dots$ *let functionals* $\phi^2_{n,n} = \phi^2_{X_1^n, X_2^n, \dots, X_n^n}$, $\delta_{n,n} = \delta_{X_1^n, X_2^n, \dots, X_n^n}$, $D^{\psi}_{n,n} = D^{\psi}_{X_1^n, X_2^n, \dots, X_n^n}$, $\rho^{(q)}_{n,n} = \rho^{(q)}_{X_1^n, X_2^n, \dots, X_n^n}$, $q \in (0, 1)$, $\mathcal{H}_{n,n} = (1/2\rho^{(q)}_{n,n})^{1/2}$, $n = 0, 1, 2, \dots$ *denote the corresponding distances. Then, as* $n \to \infty$, *if*

$$\sum_{i=1}^{n} \xi_i^n \stackrel{\mathcal{D}}{\to} Y$$

*and either* $\phi^2_{n,n} \to 0$, $\delta_{n,n} \to 0$, $D^{\psi}_{n,n} \to 0$, $\rho^{(q)}_{n,n} \to 0$ *or* $\mathcal{H}_{n,n} \to 0$ *as* $n \to \infty$, *then as* $n \to \infty$,

$$\sum_{i=1}^{n} X_i^n \stackrel{\mathcal{D}}{\to} Y.$$

*For a time series* $\{X_t\}_{t=0}^{\infty}$ *let the functionals* $\phi_t^2 = \phi^2_{X_t, X_{t+1}, \dots, X_{t+m-1}}$, $\delta_t = \delta_{X_t, X_{t+1}, \dots, X_{t+m-1}}$, $D_t^{\psi} = D^{\psi}_{X_t, X_{t+1}, \dots, X_{t+m-1}}$, $\rho_t^{(q)} = \rho^{(q)}_{X_t, X_{t+1}, \dots, X_{t+m-1}}$, $q \in (0, 1)$, $\mathcal{H}_t = (1/2\rho_t^{(q)})^{1/2}$, $t = 0, 1, 2, \dots$ *denote the* $m$*-variate Pearson coefficient, the relative entropy, the multivariate divergence measure associated with the function* $\psi$*, the generalized Tsallis entropy and the Hellinger distance for the time series, respectively.*

*Then, if, as* $t \to \infty$,

$$h(\xi_t, \xi_{t+1}, \dots, \xi_{t+m-1}) \stackrel{\mathcal{D}}{\to} Y$$

*and either* $\phi_t^2 \to 0$, $\delta_t \to 0$, $D_t^{\psi} \to 0$, $\rho_t^{(q)} \to 0$ *or* $\mathcal{H}_t \to 0$ *as* $t \to \infty$, *then, as* $t \to \infty$,

$$h(X_t, X_{t+1}, \dots, X_{t+m-1}) \stackrel{\mathcal{D}}{\to} Y.$$

From the discussion in the beginning of the present section it follows that in the case of Gaussian processes $\{X_t\}_{t=0}^{\infty}$ with $(X_t, X_{t+1}, \dots, X_{t+m-1}) \sim N(\mu_{t,m}, \Sigma_{t,m})$, the conditions of Theorem 7.2 are satisfied if, for example, $|R_{t,m}(2I_m - R_{t,m})| \to 1$ or $|\Sigma_{t,m}| / \sum_{i=0}^{m-1} \sigma^2_{t+i} \to 1$, as $t \to \infty$, where $R_{t,m}$ denote correlation matrices corresponding to $\Sigma_{t,m}$ and $(\sigma_t^2, \dots, \sigma^2_{t+m-1}) = \operatorname{diag}(\Sigma_{t,m})$. In the case of processes $\{X_t\}_{t=1}^{\infty}$ with distributions of r.v.'s $X_1, \dots, X_n$, $n \geq 1$, having generalized Eyraud–Farlie–Gumbel–Morgenstern copulas (3.3) (according to [70], this is the case for any time series of r.v.'s assuming two values), the conditions of the theorem are satisfied if, for example, $\phi_t^2 = \sum_{c=2}^{m} \sum_{i_1 < \cdots < i_c \in \{t, t+1, \dots, t+m-1\}} \alpha^2_{i_1, \dots, i_c} \to 0$ as $t \to \infty$.

It is important to emphasize here that since Theorem 7.2 holds for Tsallis entropy and multivariate divergence measures, it allows one to study convergence of statistics of time series in the case when only lower moments of the above-mentioned $U$-statistics underlying the dependence generating process for the time series exist.

Therefore, they provide a unifying approach to studying convergence in "heavy-tailed" situations and "standard" cases connected with the convergence of Pearson coefficient and the mutual information and entropy (corresponding, respectively, to the cases of second moments of the $U$-statistics and the first moments multiplied by logarithm).

The following theorem provides an estimate for the distance between the distribution function of an arbitrary statistic in dependent r.v.'s and the distribution function of the statistic in independent copies of the r.v.'s. The inequality complements (and can be better than) the well-known Pinsker's inequality for total variation between the densities of dependent and independent r.v.'s in terms of the relative entropy (see, e.g., [58]).

**Theorem 7.3.** *The following inequality holds for an arbitrary statistic* $h(X_1, \dots, X_n)$:

$$\begin{aligned} &|P(h(X_1, \dots, X_n) \le x) - P(h(\xi_1, \dots, \xi_n) \le x)| \\ &\quad \le \phi_{X_1,\dots,X_n} \max \left[ (P(h(\xi_1, \dots, \xi_n) \le x))^{1/2}, (P(h(\xi_1, \dots, \xi_n) > x))^{1/2} \right], \\ &x \in \mathbf{R}. \end{aligned}$$

The following theorems allow one to reduce the problems of evaluating expectations of general statistics in dependent r.v.'s $X_1, \dots, X_n$ to the case of independence. The theorems contain complete decoupling results for statistics in dependent r.v.'s using the relative entropy and the multivariate Pearson's $\phi^2$ coefficient. The results provide generalizations of earlier known results on complete decoupling of r.v.'s from particular dependence classes, such as martingales and adapted sequences of r.v.'s to the case of arbitrary dependence.

**Theorem 7.4.** *If* $f : \mathbf{R}^n \to \mathbf{R}$ *is a nonnegative function, then the following sharp inequalities hold:*

$$Ef(X_1, \dots, X_n) \le Ef(\xi_1, \dots, \xi_n) + \phi_{X_1,\dots,X_n} (Ef^2(\xi_1, \dots, \xi_n))^{1/2}, \tag{7.7}$$

$$Ef(X_1, \dots, X_n) \le (1 + \phi^2_{X_1,\dots,X_n})^{1/q} (Ef^q(\xi_1, \dots, \xi_n))^{1/q}, \ q \ge 2, \tag{7.8}$$

$$Ef(X_1, \dots, X_n) \le E \exp(f(\xi_1, \dots, \xi_n)) - 1 + \delta_{X_1,\dots,X_n}, \tag{7.9}$$

$$Ef(X_1, \dots, X_n) \le (1 + D^{\psi}_{X_1,\dots,X_n})^{(1-\frac{1}{q})} (Ef^q(\xi_1, \dots, \xi_n))^{1/q}, \ q > 1, \tag{7.10}$$

*where* $\psi(x) = |x|^{q/(q-1)} - 1$.

**Remark 7.1.** It is interesting to note that from relation (7.2) and inequality (7.7) it follows that the following representation holds for the multivariate Pearson coefficient $\phi_{X_1,\dots,X_n}$:

$$\phi_{X_1,\dots,X_n} = \max_{\substack{f: Ef(\xi_1,\dots,\xi_n)=0, \\ Ef^2(\xi_1,\dots,\xi_n)<\infty}} \frac{(Ef(X_1, \dots, X_n) - Ef(\xi_1, \dots, \xi_n))}{(Ef^2(\xi_1, \dots, \xi_n))^{1/2}}. \tag{7.11}$$

The following result gives complete decoupling inequalities for the tail probabilities of arbitrary statistics in dependent r.v.'s.

**Theorem 7.5.** *The following inequalities hold:*

$$P\left(h(X_1,\ldots,X_n)>x\right)\le P(h(\xi_1,\ldots,\xi_n)>x)+\phi_{X_1,\ldots,X_n}\left(P(h(\xi_1,\ldots,\xi_n)>x)\right)^{\frac{1}{2}},$$

$$P\left(h(X_1,\ldots,X_n)>x\right)\le\left(1+\phi^2_{X_1,\ldots,X_n}\right)^{1/2}\left(P(h(\xi_1,\ldots,\xi_n)>x)\right)^{1/2},$$

$$P(h(X_1,\ldots,X_n)>x)\le(e-1)P(h(\xi_1,\ldots,\xi_n)>x)+\delta_{X_1,\ldots,X_n},$$

$$P(h(X_1,\ldots,X_n)>x)\le\left(1+D^{\psi}_{X_1,\ldots,X_n}\right)^{(1-\frac{1}{q})}\left(P(h(\xi_1,\ldots,\xi_n)>x)\right)^{\frac{1}{q}},q>1,$$

$x\in\mathbf{R}$, *where* $\psi(x)=|x|^{q/(q-1)}-1$.

## 8. Appendix: Proofs

*Proof of Theorem 2.1.* Let us first prove the necessity part of the theorem. Denote

$$T(x_1,\ldots,x_n)=\int_{-\infty}^{x_1}\cdots\int_{-\infty}^{x_n}(1+U_n(t_1,\ldots,t_n))\prod_{i=1}^n dF_i(t_i).$$

Let $k\in\{1,\ldots,n\}$, $x_k\in\mathbf{R}$. Let us show that

$$T(\infty,\ldots,\infty,x_k,\infty,\ldots,\infty)=F_k(x_k), \tag{8.1}$$

$x_k\in\mathbf{R}$, $k=1,\ldots,n$. It suffices to consider the case $k=1$. We have

$$\begin{aligned}&T(x_1,\infty,\ldots,\infty)\\&\quad=\underbrace{\int_{-\infty}^{x_1}\int_{-\infty}^{\infty}\cdots\int_{-\infty}^{\infty}}_{n}(1+U_n(t_1,\ldots,t_n))\prod_{i=1}^n dF_i(t_i)\\&\quad=F_1(x_1)+\sum_{c=2}^n\sum_{1\le i_1<\cdots<i_c\le n}\underbrace{\int_{-\infty}^{x_1}\int_{-\infty}^{\infty}\cdots\int_{-\infty}^{\infty}}_{n}g_{i_1,\ldots,i_c}(t_{i_1},\ldots,t_{i_c})\prod_{i=1}^n dF_i(t_i)\\&\quad=F_1(x_1)+\Sigma''.\end{aligned}$$

It is easy to see that there is at least one $t_s$ of $t_2,\ldots,t_n$ among the arguments of each of the functions $g_{i_1,\ldots,i_c}(t_{i_1},\ldots,t_{i_c})$, $1\le i_1<\cdots<i_c\le n$, $c=2,\ldots,n$, in the latter summand. By A2 we get, therefore, that $\Sigma''=0$. Consequently, $T(x_1,\infty,\ldots,\infty)=F_1(x_1)$, $x_1\in\mathbf{R}$, and (8.1) holds. It is evident that

$$\lim_{x_k\to-\infty}T(x_1,\ldots,x_k,\ldots,x_n)=0 \tag{8.2}$$

for all $x_j\in\mathbf{R}$, $j=1,\ldots,n$, $j\ne k$, $k=1,\ldots,n$. Since

$$T(x_1,\ldots,x_n)=\prod_{i=1}^n F_i(x_i)+E\left[U_n(\xi_1,\ldots,\xi_n))\prod_{i=1}^n I(\xi_i\le x_i)\right],$$

from the monotone convergence theorem we obtain that $T(x_1,\ldots,x_n)$ is right-continuous in $(x_1,\ldots,x_n)\in\mathbf{R}^n$. Let $\delta^k_{[a,b)}T(x_1,\ldots,x_n)=T(x_1,\ldots,x_{k-1},b,x_{k+1},$

$\ldots, x_n) - T(x_1, \ldots, x_{k-1}, a, x_{k+1}, \ldots, x_n)$, $a < b$. By integrability of the functions $g_{i_1,\ldots,i_c}$ and condition A3 we obtain ($I(\cdot)$ denotes the indicator function)

$$\begin{aligned}(8.3)\quad & \delta^1_{(a_1,b_1]}\delta^2_{(a_2,b_2]}\cdots\delta^n_{(a_n,b_n]}T(x_1,\ldots,x_n)\\ &= \prod_{i=1}^n P(a_i < \xi_i \le b_i) + E\left[U_n(\xi_{i_1},\ldots,\xi_{i_n})\prod_{i=1}^n I(a_i < \xi_i \le b_i)\right] \ge 0\end{aligned}$$

for all $a_i < b_i$, $i = 1, \ldots, n$.[1] Right-continuity of $T(x_1, \ldots, x_n)$ and (8.1)–(8.3) imply that $T(x_1, \ldots, x_n)$ is a joint cdf of some r.v.'s $X_1, \ldots, X_n$ with one-dimensional cdf's $F_k(x_k)$, and the joint cdf $T(x_1, \ldots, x_n)$ satisfies (2.1).

Let us now prove the sufficiency part. Consider the functions

$$f_{i_1,\ldots,i_c}(x_{i_1},\ldots,x_{i_c}) = \sum_{s=2}^c (-1)^{c-s} \sum_{j_1<\cdots<j_s\in\{i_1,\ldots,i_c\}} \left(\frac{dF(x_{j_1},\ldots,x_{j_s})}{dF_{j_1}\cdots dF_{j_s}} - 1\right),$$

$1 \le i_1 < \cdots < i_c \le n$, $c = 2, \ldots, n$. Obviously, the functions $f_{i_1,\ldots,i_c}$ satisfy condition A1. Let us show that they satisfy condition A2. It suffices to consider the case $i_1 = 1, i_2 = 2, \ldots, i_c = c, k = 1$. We have

$$\begin{aligned}&Eg_{1,2,\ldots,c}(\xi_1, x_2, \ldots, x_c)\\ &\quad= \int_{-\infty}^{\infty} g_{1,2,\ldots,c}(x_1, x_2, \ldots, x_c)dF_1(x_1)\\ &\quad= \int_{-\infty}^{\infty}\sum_{s=2}^c(-1)^{c-s}\left[\sum_{2\le i_2<\cdots<i_s\le c}\left(\frac{dF(x_1, x_{i_2},\ldots,x_{i_s})}{dF_1 dF_{i_2}\cdots dF_{i_s}} - 1\right)\right.\\ &\qquad\qquad\left. + \sum_{2\le i_1<\cdots<i_s\le c}\left(\frac{dF(x_{i_1}, x_{i_2},\ldots,x_{i_s})}{dF_{i_1} dF_{i_2}\cdots dF_{i_s}} - 1\right)\right] dF_1(x_1)\\ &\quad= \sum_{s=2}^c(-1)^{c-s}\left\{\sum_{2\le i_2<\cdots<i_s\le c}(\frac{dF(x_{i_2},\ldots,x_{i_s})}{dF_{i_2}\cdots dF_{i_s}} - 1)\right.\\ &\qquad\qquad\left. + \sum_{2\le i_1<\cdots<i_s\le c}(\frac{dF(x_{i_1},\ldots,x_{i_s})}{dF_{i_1}\cdots dF_{i_s}} - 1)\right\} = 0.\end{aligned}$$

By the inversion formula (see, e.g., [8, pp. 177–178]) it follows that if $a_{i_1,\ldots,i_c}$, $b_{i_1,\ldots,i_c}$, $1 \le i_1 < \cdots < i_c \le n$, $c = 2, \ldots, n$, are arbitrary numbers then the relations

$$b_{i_1,\ldots,i_c} = \sum_{c=2}^n \sum_{j_1<\cdots<j_s\in\{i_1,\ldots,i_c\}} a_{j_1,\ldots,j_s}, 1 \le i_1 < \cdots < i_c \le n, c = 2, \ldots, n,$$

and

$$a_{i_1,\ldots,i_c} = \sum_{s=2}^c (-1)^{c-s} \sum_{j_1<\cdots<j_s\in\{i_1,\ldots,i_c\}} b_{i_1,\ldots,i_c}, 1 \le i_1 < \cdots < i_c \le n, c = 2, \ldots, n,$$

---

[1] Note that (8.2) and (8.3) are immediate if the probability space and the random variables are defined in the canonical way with $\Omega = \mathbf{R}^n$ and $X_i(\omega) = \omega_i$ for $\omega = (\omega_1, \ldots, \omega_n)$ and $P(A) = \int_A (1 + \sum_{1\le i_1<\cdots<i_c\le n} g_{i_1,\ldots,i_c}(t_{i_1}, \ldots, t_{i_c})) \prod_{i=1}^n dF_i(t_i)$.

are equivalent. Taking $a_{i_1,\dots,i_c} = g_{i_1,\dots,i_c}(x_{i_1},\dots,x_{i_c})$, $b_{i_1,\dots,i_c} = dF(x_{i_1},\dots,x_{i_c})/\prod_{j=1}^c dF_{i_j} - 1$, for $1 \le i_1 < \cdots < i_c \le n, c = 2,\dots,n$, we obtain, in particular, that

$$G(x_1,\dots,x_n) = 1 + U_n(x_1,\dots,x_n). \tag{8.4}$$

Therefore, representation (2.1) holds and the functions $f_{i_1,\dots,i_c}$ satisfy condition A3. Suppose now that there exists another set of functions $g_{i_1,\dots,i_c}$ satisfying conditions A1–A3 and such that (2.1) holds and, equivalently, (2.2) holds. Let $B_s$ be a set of $s$ integers $j_1,\dots,j_s$ with $1 \le j_1 < \cdots < j_s \le n$, for $s = 2,\dots,n$. By Remark 1 we have

$$\sum_{c=2}^{s} \sum_{\{i_1<\cdots<i_c\}\in B_s} f_{i_1,\dots,i_c}(\xi_{i_1},\dots,\xi_{i_c}) = \sum_{c=2}^{s} \sum_{\{i_1<\cdots<i_c\}\in B_s} g_{i_1,\dots,i_c}(\xi_{i_1},\dots,\xi_{i_c}) \tag{8.5}$$

(a.s.). From (8.5) we subsequently obtain that $f_{i_1,i_2}(\xi_{i_1},\xi_{i_2}) = g_{i_1,i_2}(\xi_{i_1},\xi_{i_2})$ (a.s.), $1 \le i_1 < i_2 \le n$; $f_{i_1,i_2,i_3}(\xi_{i_1},\xi_{i_2},\xi_{i_3}) = g_{i_1,i_2,i_3}(\xi_{i_1},\xi_{i_2},\xi_{i_3})$ (a.s.), $1 \le i_1 < i_2 < i_3 \le n$; ..., $f_{1,2,\dots,n}(\xi_1,\xi_2,\dots,\xi_n) = g_{1,2,\dots,n}(\xi_1,\xi_2,\dots,\xi_n)$ (a.s.), that is $g_{i_1,\dots,i_c}(\xi_{i_1},\dots,\xi_{i_c}) = f_{i_1,\dots,i_c}(\xi_{i_1},\dots,\xi_{i_c})$ (a.s.), $1 \le i_1 < \cdots < i_c \le n$, $c = 2,\dots,n$. This completes the proof. □

*Proof of Theorems 3.2–4.1.* By definition, a function $C : [0,1]^n \to [0,1]$ is a $n$-dimensional copula if and only if it is a joint cdf of $n$ r.v.'s each of which is uniformly distributed on $[0,1]$. Let, as in Section 3, $V_1,\dots,V_n$ denote independent r.v.'s with uniform distribution on $[0,1]$. From Theorem 2.1 we obtain that $C : [0,1]^n \to [0,1]$ is an absolutely continuous copula if and only if there exist functions $\tilde{g}_{i_1,\dots,i_c} : \mathbf{R}^c \to \mathbf{R}$, $1 \le i_1 < \cdots < i_c \le n$, $c = 2,\dots,n$, satisfying the conditions

$$E|\tilde{g}_{i_1,\dots,i_c}(V_{i_1},\dots,V_{i_c})| < \infty,$$

$$E(\tilde{g}_{i_1,\dots,i_c}(V_{i_1},\dots,V_{i_{k-1}},V_{i_k},V_{i_{k+1}},\dots,V_{i_c})|V_{i_1},\dots,V_{i_{k-1}},V_{i_{k+1}},\dots,V_{i_c}) = 0\ (a.s.),$$

$1 \le i_1 < \cdots < i_c \le n$, $k = 1,2,\dots,c$, $c = 2,\dots,n$,

$$\sum_{c=2}^{n} \sum_{1\le i_1<\dots<i_c\le n} \tilde{g}_{i_1,\dots,i_c}(V_{i_1},\dots,V_{i_c}) \ge -1\ (a.s.),$$

and such that representation (3.1) holds. This proves Theorem 3.2. Theorem 3.3 follows from Theorem 3.2 and the relation (given by Sklar's Theorem 3.1) $F_{X_1,\dots,X_n}(x_1,\dots,x_n) = C_{X_1,\dots,X_n}(F_{X_1}(x_1),\dots,F_{X_n}(x_n))$, $x_k \in \mathbf{R}$, $k = 1,\dots,n$, between the joint distribution functions $F_{X_1,\dots,X_n}(x_1,\dots,x_n)$ and the corresponding copulas $C_{X_1,\dots,X_n}(u_1,\dots,u_n)$. If $U_n$ is the $U$-statistic corresponding to the r.v.'s $X_1,\dots,X_n$, then for any Borel measurable function $f$ for which the expectations exist, one has

$$\begin{aligned} Ef(X_1,\dots,X_n) &= E\{f(\xi_1,\dots,\xi_n)(1 + \sum_{c=2}^{n} \sum_{1\le i_1<\cdots<i_c\le n} g_{i_1,\dots,i_c}(\xi_{i_1},\dots,\xi_{i_c}))\} \\ &= E\left[f(\xi_1,\dots,\xi_n)(1 + U_n(\xi_1,\dots,\xi_n))\right]. \end{aligned}$$

This and Theorem 2.1 implies Theorem 4.1. □

*Proof of Theorem 5.1.* It is evident that $g_{i_1,\dots,i_c}(x_{i_1},\dots,x_{i_c}) = 0$ satisfy conditions A1–A3. If $g_{i_1,\dots,i_c}(x_{i_1},\dots,x_{i_c}) = 0$, $1 \le i_1 < \cdots < i_c \le n$, $c = 2,\dots,n$, then representation (2.1) takes the form

$$F(x_1,\dots,x_n) = \prod_{i=1}^{n} F_i(x_i). \tag{8.6}$$

R.v's with the joint distribution function (8.6) are independent. Let now $X_1, \dots, X_n$ be independent r.v.'s with one-dimensional distribution functions $F_i(x_i), i = 1, \dots, n$. Then their joint distribution function has form (8.6). This and the uniqueness of the functions $g_{i_1,\dots,i_c}$ given by Theorem 2.1 completes the proof of the theorem. □

*Proof of Theorems 5.2–5.8.* Below, we give proofs of Theorems 5.7 and 5.8. The rest of the theorems can be proven in a similar way. Let $X_1, \dots, X_n$ be r.v.'s with the joint distribution function satisfying representation (2.2) with functions $g_{i_1,\dots,i_c}$ such that $E\xi_{i_1}^{\alpha_{i_1}} \cdots \xi_{i_c}^{\alpha_{i_c}} g_{i_1,\dots,i_c}(\xi_{i_1}, \dots, \xi_{i_c}) = 0$, $1 \le i_1 < \cdots < i_c \le n$, $c = 2, \dots, n$, where $\xi_1, \dots, \xi_n$ are independent copies of $X_k$, $k = 1, 2, \dots, n$. Let $\alpha_j \in \{0, 1, \dots, \alpha\}$, $j = 1, \dots, n$. Taking in (4.2) $f(x_1, \dots, x_n) = \prod_{j=1}^n x_j^{\alpha_j}$ and using independence of the r.v.'s $\xi_1, \dots, \xi_n$ we obtain

$$\begin{aligned}
E\prod_{j=1}^n X_j^{\alpha_j} &= E\prod_{j=1}^n \xi_j^{\alpha_j} + \sum_{c=2}^n \sum_{1\le i_1<\cdots<i_c\le n} E\prod_{j=1}^n \xi_j^{\alpha_j} g_{i_1,\dots,i_c}(\xi_{i_1}, \dots, \xi_{i_c}) \\
&= \prod_{j=1}^n E\xi_j^{\alpha_j} + \sum_{c=2}^n \sum_{1\le i_1<\cdots<i_c\le n} E\prod_{k=1}^c \xi_{i_k}^{\alpha_{i_k}} g_{i_1,\dots,i_c}(\xi_{i_1}, \dots, \xi_{i_c}) \\
&\qquad \times \prod_{\substack{k=1,\dots,n\\ k\ne i_1,\dots,i_c}} E\xi_k^{\alpha_k}.
\end{aligned} \tag{8.7}$$

Using $E\xi_{i_1}^{\alpha_{i_1}} \cdots \xi_{i_c}^{\alpha_{i_c}} g_{i_1,\dots,i_c}(\xi_{i_1}, \dots, \xi_{i_c}) = 0$, we get from (8.7) $E\prod_{j=1}^n X_j^{\alpha_j} = \prod_{j=1}^n E\xi_j^{\alpha_j} = \prod_{j=1}^n EX_j^{\alpha_j}$, that is the r.v.'s $X_1, \dots, X_n$ form a multiplicative system of order $\alpha$. Let us suppose now that r.v.'s $X_1, \dots, X_n$ form a multiplicative system of order $\alpha$, that is, $E|X_j|^\alpha < \infty$, $j = 1, \dots, n$, and for all $\alpha_j \in \{0, 1, \dots, \alpha\}$, $j = 1, \dots, n$, $E\prod_{j=1}^n X_j^{\alpha_j} = \prod_{j=1}^n EX_j^{\alpha_j}$. From Remark 2.1 and Theorem 4.1 it follows that

$$E\prod_{r=1}^k X_{j_r}^{\alpha_{j_r}} = \prod_{r=1}^k EX_{j_r}^{\alpha_{j_r}} + \sum_{c=2}^k \sum_{i_1<\cdots<i_c\in\{j_1,\dots,j_k\}} \prod_{r=1}^k \xi_{j_r}^{\alpha_{j_r}} g_{i_1,\dots,i_c}(\xi_{i_1}, \dots, \xi_{i_c}),$$

$\alpha_{j_r} \in \{0, 1, \dots, \alpha\}$, $1 \le j_1 < \cdots < j_k \le n$, $k = 2, \dots, n$. Therefore, for all $\alpha_{j_r} \in \{0, 1, \dots, \alpha\}$, $1 \le j_1 < \cdots < j_k \le n$, $k = 2, \dots, n$,

$$\sum_{c=2}^k \sum_{i_1<\cdots<i_c\in\{j_1,\dots,j_k\}} E\prod_{r=1}^c \xi_{i_r}^{\alpha_{i_r}} g_{i_1,\dots,i_c} \prod_{\substack{r=1,\dots,k\\ j_r\ne i_1,\dots,i_c}} E\xi_{j_r}^{\alpha_{j_r}} = 0. \tag{8.8}$$

From (8.8) we subsequently obtain that

$$E\xi_{i_1}^{\alpha_{i_1}}\xi_{i_2}^{\alpha_{i_2}} g_{i_1,i_2}(\xi_{i_1}, \xi_{i_2}) = 0, \alpha_{i_1}, \alpha_{i_2} \in \{0, 1, \dots, \alpha\}, 1 \le i_1 < i_2 \le n;$$

$$E\xi_{i_1}^{\alpha_{i_1}}\xi_{i_2}^{\alpha_{i_2}}\xi_{i_3}^{\alpha_{i_3}} g_{i_1,i_2,i_3}(\xi_{i_1}, \xi_{i_2}, \xi_{i_3}) = 0, \alpha_{i_1}, \alpha_{i_2}, \alpha_{i_3} \in \{0, 1, \dots, \alpha\}, 1 \le i_1 < i_2$$
$$< i_3 \le n \text{ and } E\xi_1^{\alpha_1} \cdots \xi_n^{\alpha_n} g_{1,2,\dots,n}(\xi_1, \xi_2, \dots, \xi_n), \alpha_k \in \{0, 1, \dots, \alpha\}, k = 1, \dots, n.$$

Therefore, $E\xi_{i_1}^{\alpha_{i_1}} \cdots \xi_{i_c}^{\alpha_{i_c}} g_{i_1,\dots,i_c}(\xi_{i_1}, \dots, \xi_{i_c}) = 0$, $1 \le i_1 < \cdots < i_c \le n$, $c = 2, \dots, n$. Let $X_1, \dots, X_n$ be $r$-independent r.v.'s. From Remark 2.1 it follows that

$$\begin{aligned}
&F(x_{j_1}, \dots, x_{j_k}) \\
&\quad = \int_{-\infty}^{x_{j_1}} \cdots \int_{-\infty}^{x_{j_k}} \Big(1 + \sum_{c=2}^n \sum_{i_1<\ldots<i_c\in B_k} g_{i_1,\dots,i_c}(t_{i_1}, \dots, t_{i_c})\Big) \prod_{i=1}^k dF_{j_i}(t_{j_i})
\end{aligned}$$

for all $B_k = \{1 \le j_1 < \cdots < j_k \le n\}$, $k = 1, \dots, r$. Using $r$-independence of $X_i'$s we subsequently obtain from here that $g_{i_1,i_2}(\xi_{i_1}, \xi_{i_2}) = 0$ (a.s.), $1 \le i_1 < i_2 \le n$; $g_{i_1,i_2,i_3}(\xi_{i_1}, \xi_{i_2}, \xi_{i_3}) = 0$ (a.s.), $1 \le i_1 < i_2 < i_3 \le n; \dots, g_{i_1,\dots,i_r}(\xi_{i_1}, \dots, \xi_{i_r}) = 0$ (a.s.), $1 \le i_1 < \cdots < i_r \le n$. Therefore, $g_{i_1,\dots,i_c}(\xi_{i_1}, \dots, \xi_{i_c}) = 0$ (a.s.), $1 \le i_1 < \cdots < i_c \le n$, $c = 2, \dots, r$. Let now $X_1, \dots, X_n$ be r.v.'s such that the functions $g_{i_1,\dots,i_c}$ in representations (2.1) and (2.2) satisfy the conditions $g_{i_1,\dots,i_c}(\xi_{i_1}, \dots, \xi_{i_c}) = 0$ (a.s.), $1 \le i_1 < \cdots < i_c \le n$, $c = 2, \dots, r$, that is their joint distribution function has the form

$$\begin{aligned} &F(x_1, \dots, x_n) \\ &\quad = P(X_1 \le x_1, \dots, X_n \le x_n) \\ &\quad = \int_{-\infty}^{x_1} \cdots \int_{-\infty}^{x_n} (1 + \sum_{c=r+1}^{n} \sum_{1 \le i_1 < \cdots < i_c \le n} g_{i_1,\dots,i_c}(t_{i_1}, \dots, t_{i_c})) \prod_{i=1}^{n} F_i(t_i). \end{aligned}$$

Let $1 \le j_1 < \cdots < j_r \le n$. Let us show that the r.v.'s $X_{j_1}, \dots, X_{j_r}$ are jointly independent. Without loss of generality, it suffices to consider the case $j_1 = 1, \dots, j_r = r$. We have, similar to the proof of Theorem 2.1,

$$\begin{aligned} &F(x_1, \dots, x_r) \\ &\quad = P(X_1 \le x_1, \dots, X_r \le x_r) \\ &\quad = \underbrace{\int_{-\infty}^{x_1} \cdots \int_{-\infty}^{x_r} \int_{-\infty}^{\infty} \int_{-\infty}^{\infty}}_{n} \left(1 + \sum_{c=r+1}^{n} \sum_{1 \le i_1 < \cdots < i_c \le n} g_{i_1,\dots,i_c}(t_{i_1}, \dots, t_{i_c})\right) \\ &\qquad\qquad\qquad \times \prod_{i=1}^{n} dF_i(t_i) \\ &\quad = \prod_{i=1}^{r} F_i(x_i) + \sum_{c=r+1}^{n} \sum_{1 \le i_1 < \cdots < i_c \le n} g_{i_1,\dots,i_c}(t_{i_1}, \dots, t_{i_c}) \prod_{i=1}^{n} F_i(t_i) \\ &\quad = \prod_{i=1}^{r} F_i(x_i) + \Sigma''. \end{aligned}$$

It is easy to see that there is at least one $t_s$ of $t_{r+1}, \dots, t_n$ among the arguments of each of the functions $g_{i_1,\dots,i_c}(t_{i_1}, \dots, t_{i_c})$, $1 \le i_1 < \cdots < i_c \le n$, $c = r+1, \dots, n$, in the latter summand and, therefore, by A2, $\Sigma'' = 0$. Consequently, $F(x_1, \dots, x_r) = \prod_{i=1}^{r} F_i(x_i)$ that establishes joint independence of the r.v.'s $X_1, \dots, X_r$. The proof is complete. □

*Proof of Theorem 6.1.* Evidently, if the r.v.'s $X_1, \dots, X_n$ are jointly independent, then they form a multiplicative system of order $\alpha$. Let us show that if $card(A_i) \le \alpha + 1$, $i = 1, \dots, n$, and r.v.'s $X_1, \dots, X_n$ form a multiplicative system of order $\alpha$, then they are jointly independent. It suffices to show that

$$E \prod_{i=1}^{n} f_i(X_i) = \prod_{i=1}^{n} E f_i(X_i) \tag{8.9}$$

for all continuous functions $f_i : \mathbf{R} \to \mathbf{R}$, $i = 1, \dots, n$, vanishing outside a finite interval. Let $1 \le i \le n$. It is easy to see that if $f_i(x)$, $x \in \mathbf{R}$, is an arbitrary function, then there exists a polynomial $r_i(x)$, $x \in \mathbf{R}$, of degree not greater than $\alpha$, such that $f_i(x) = r_i(x)$, $x \in A_i$, that is $f_i(\xi_i) = r_i(\xi_i)$ (a.s.). Using Theorems 4.1 and 5.7, we get that for all continuous functions $f_i : \mathbf{R} \to \mathbf{R}$ (below, $r_i(x)$ are polynomials corresponding to $f_i(x)$)

$$\begin{aligned} E\prod_{i=1}^{n} f_i(X_i) &= E\prod_{i=1}^{n} f_i(\xi_i) + \sum_{c=2}^{n} \sum_{1\le i_1<\cdots<i_c\le n} E\prod_{i=1}^{n} f_i(\xi_i) g_{i_1,\dots,i_c}(\xi_{i_1},\dots,\xi_{i_c}) \\ &= E\prod_{i=1}^{n} f_i(\xi_i) + \sum_{c=2}^{n} \sum_{1\le i_1<\cdots<i_c\le n} E\prod_{i=1}^{n} r_i(\xi_i) g_{i_1,\dots,i_c}(\xi_{i_1},\dots,\xi_{i_c}) \\ &= E\prod_{i=1}^{n} f_i(\xi_i). \end{aligned}$$

The proof is complete. □

*Proof of Theorem 7.1.* From concavity of the functions $\log(1+x)$, relation (7.2) and the inequality $\log(1+x) \le x$, $x \ge 0$, we get $\delta_{X_1,\dots,X_n} = E\log(1 + U_n(X_1,\dots,X_n)) \le \log(1 + EU_n(X_1,\dots,X_n)) = \log(1+\phi^2_{X_1,\dots,X_n}) \le \phi^2_{X_1,\dots,X_n}$. The proof is complete. □

*Proof of Theorem 7.2.* Let $h(\xi_t,\xi_{t+1},\dots,\xi_{t+m-1}) \xrightarrow{\mathcal{D}} Y$. By Theorem 4.1 we have that for any continuous bounded function $g : \mathbf{R} \to \mathbf{R}$

$$Eg(h(X_t,\dots,X_{t+m-1})) = Eg(h(\xi_t,\dots,\xi_{t+m-1}))(1 + U_m(\xi_t,\dots,\xi_{t+m-1})).$$

By Chebyshev's inequality and (7.2) we have that for all $\epsilon > 0$,

$$P(|U_m(\xi_t,\xi_{t+1},\dots,\xi_{t+m-1})| > \epsilon) \le \phi_t^2/\epsilon^2.$$

Since the function $w(x) = (1+x)\ln(1+x) - x$ is increasing in $x \in [0,\infty)$ and decreasing in $x \in (-\infty, 0)$, we have that if $U_m(\xi_t,\xi_{t+1},\dots,\xi_{t+m-1}) > \epsilon$ or $U_m(\xi_t,\xi_{t+1},\dots,\xi_{t+m-1}) < -\epsilon$, then

$$w(U_m(\xi_t,\xi_{t+1},\dots,\xi_{t+m-1})) > (w(\epsilon) \wedge w(-\epsilon)),$$

where $a \wedge b = \min(a,b)$. Therefore, by Chebyshev's inequality, (7.1) and since $EU_m(\xi_t,\dots,\xi_{t+m-1}) = 0$ (by condition A2) we get, for $0 < \epsilon < 1$,

$$\begin{aligned} &P(|U_m(\xi_t,\dots,\xi_{t+m-1})| > \epsilon) \\ &\quad \le P(w(U_m(\xi_t,\dots,\xi_{t+m-1})) > (w(\epsilon) \wedge w(-\epsilon))) \\ &\quad \le Ew(U_m(\xi_t,\xi_{t+1},\dots,\xi_{t+m-1}))/(w(\epsilon) \wedge w(-\epsilon)) \\ &\quad = E\,(1 + U_m(\xi_t,\dots,\xi_{t+m-1})) \\ &\qquad \times \log(1 + U_m(\xi_t,\dots,\xi_{t+m-1}))\,/(w(\epsilon) \wedge w(-\epsilon)) \\ &\quad = \delta_t/(w(\epsilon) \wedge w(-\epsilon)). \end{aligned} \tag{8.10}$$

If $\epsilon \ge 1$, Chebyshev's inequality and $U_m(\xi_t,\dots,\xi_{t+m-1}) \ge -1$ yield

$$P(|U_m(\xi_t,\dots,\xi_{t+m-1})| > \epsilon) \le \frac{Ew(U_m(\xi_t,\dots,\xi_{t+m-1}))}{w(\epsilon)} = \delta_t/w(\epsilon). \tag{8.11}$$

Similar to the above, by Chebyshev's inequality and (7.3), for $0 < \epsilon < 1$,

$$\begin{aligned} P(|U_m(\xi_t,\dots,\xi_{t+m-1})| > \epsilon) &\le P(\psi(1+U_m(\xi_t,\dots,\xi_{t+m-1})) > (\psi(1+\epsilon) \wedge \psi(1-\epsilon)) \\ &\le E\psi(1 + U_m(\xi_t,\dots,\xi_{t+m-1}))/(\psi(1+\epsilon) \wedge \psi(1-\epsilon)) \\ &= D_t^{\psi}/(\psi(1+\epsilon) \wedge \psi(1-\epsilon)). \end{aligned} \tag{8.12}$$

For $\epsilon \geq 1$,

$$P(|U_m(\xi_t, \ldots, \xi_{t+m-1})| > \epsilon) \leq P(\psi(1 + U_m(\xi_t, \ldots, \xi_{t+m-1})) > \psi(1 + \epsilon)) \leq D_t^{\psi}/\psi(1 + \epsilon). \tag{8.13}$$

Inequalities (8.10)–(8.13) imply that $U_m(\xi_t, \xi_{t+1}, \ldots, \xi_{t+m-1}) \to 0$ (in probability) as $t \to \infty$, if $\phi_t^2 \to 0$, or $\delta_t \to 0$, or $D_t^{\psi} \to 0$ as $t \to \infty$. The same argument as in the case of the measure $D_t^{\psi}$, used with $\psi(x) = x^{1-q}$, establishes that $U_m(\xi_t, \xi_{t+1}, \ldots, \xi_{t+m-1}) \to 0$ (in probability) as $t \to \infty$, if $\rho_t^{(q)} \to 0$ as $t \to \infty$ for $q \in (0, 1)$. In particular, the latter holds for the case $q = 1/2$, and, consequently, for the Hellinger distance $\mathcal{H}_t$. The above implies, by Slutsky theorem, that $Eg(h(X_t, X_{t+1}, \ldots, X_{t+m-1})) \to Eg(Y)$ as $t \to \infty$. Since this holds for any continuous bounded function $g$, we get $h(X_t, X_{t+1}, \ldots, X_{t+m-1}) \to Y$ (in distribution) as $t \to \infty$. The proof is complete. The case of double arrays requires only minor notational modifications. □

*Proof of Theorem 7.3.* From Theorem 4.1, relation (7.2) and Hölder inequality we obtain that for any $x \in \mathbf{R}$ and r.v.'s $X_1, \ldots, X_n$

$$\begin{aligned} &P(h(X_1, \ldots, X_n) \leq x) - P(h(\xi_1, \ldots, \xi_n) \leq x) \\ &\quad = EI(h(\xi_1, \ldots, \xi_n) \leq x) U_n(\xi_1, \ldots, \xi_n) \\ &\quad \leq \phi_{X_1, \ldots X_n}(P(h(\xi_1, \ldots, \xi_n) \leq x))^{1/2}, \end{aligned} \tag{8.14}$$

$$\begin{aligned} &P(h(X_1, \ldots, X_n) > x) - P(h(\xi_1, \ldots, \xi_n) > x) \\ &\quad = EI(h(\xi_1, \ldots, \xi_n) > x) U_n(\xi_1, \ldots, \xi_n) \\ &\quad \leq \phi_{X_1, \ldots, X_n}(P(h(\xi_1, \ldots, \xi_n) > x))^{1/2}. \end{aligned} \tag{8.15}$$

The latter inequalities imply that for any $x \in \mathbf{R}$

$$\begin{aligned} &|P(h(X_1, \ldots, X_n) \leq x) - P(h(\xi_1, \ldots, \xi_n) \leq x)| \\ &\quad \leq \phi_{X_1, \ldots, X_n} \max \left[ (P(h(\xi_1, \ldots, \xi_n) \leq x))^{1/2}, (P(h(\xi_1, \ldots, \xi_n) > x))^{1/2} \right]. \end{aligned}$$

The proof is complete. □

*Proof of Theorem 7.4.* By Theorem 4.1 we have $Ef(X_1, \ldots, X_n) = Ef(\xi_1, \ldots, \xi_n) + EU_n(\xi_1, \ldots, \xi_n) f(\xi_1, \ldots, \xi_n)$. By Cauchy–Schwarz inequality and relation (7.2) we get

$$EU_n(\xi_1, \ldots, \xi_n) f(\xi_1, \ldots, \xi_n) \leq \left( EU_n^2(\xi_1, \ldots, \xi_n) \right)^{1/2} \left( Ef^2(\xi_1, \ldots, \xi_n) \right)^{1/2}.$$

Therefore, (7.7) holds. Sharpness of (7.7) follows from the choice of independent $X_1, \ldots, X_n$. Similarly, from Hölder inequality it follows that if $q > 1$, $1/p + 1/q = 1$, then

$$Ef(X_1, \ldots, X_n) \leq (E(1 + U_n(\xi_1, \ldots, \xi_n))^p)^{1/p} (Ef(\xi_1, \ldots, \xi_n))^q)^{1/q}. \tag{8.16}$$

This implies (7.10). If in estimate (8.16) $q \geq 2$ and, therefore, $p \in (1, 2]$, by Theorem 4.1, Jensen inequality and relation (7.2) we have

$$\begin{aligned} E(1 + U_n(\xi_1, \ldots, \xi_n))^p &= E(1 + U_n(X_1, \ldots, X_n))^{p-1} \\ &\leq (1 + EU_n(X_1, \ldots, X_n))^{p-1} \\ &= (1 + \phi_{X_1, \ldots, X_n}^2)^{p/q}. \end{aligned}$$

Therefore, (7.8) holds. Sharpness of (7.8) and (7.10) follows from the choice of $X_i = const$ (a.s.), $i = 1, \ldots, n$. According to Young's inequality (see [19, p. 512]), if $p : [0, \infty) \to [0, \infty)$ is a non-decreasing right-continuous function satisfying $p(0) = \lim_{t \to 0+} p(t) = 0$ and $p(\infty) = \lim_{t \to \infty} p(t) = \infty$, and $q(t) = \sup\{u : p(u) \leq t\}$ is a right-continuous inverse of $p$, then

$$st \leq \phi(s) + \psi(t), \tag{8.17}$$

where $\phi(t) = \int_0^t p(s)ds$ and $\psi(t) = \int_0^t q(s)ds$. Using (8.17) with $p(t) = ln(1+t)$ and (7.1), we get that

$$\begin{aligned} EU_n(\xi_1, \ldots, \xi_n) f(\xi_1, \ldots, \xi_n) &\leq E(e^{f(\xi_1, \ldots, \xi_n)}) - 1 - Ef(\xi_1, \ldots, \xi_n) \\ &\quad + E(1 + U_n(\xi_1, \ldots, \xi_n)) \log(1 + U_n(\xi_1, \ldots, \xi_n)) \\ &= E(e^{f(\xi_1, \ldots, \xi_n)}) - 1 - Ef(\xi_1, \ldots, \xi_n) \\ &\quad + \delta_{X_1, \ldots, X_n}. \end{aligned}$$

This establishes (7.9). Sharpness of (7.9) follows, e.g., from the choice of independent $X_i'$s and $f \equiv 0$. □

*Proof of Theorem 7.5.* The theorem follows from inequalities (7.7)–(7.10) applied to $f(x_1, \ldots, x_n) = I(h(x_1, \ldots, x_n) > x)$. □

## Acknowledgements

The authors are grateful to Peter Phillips, two anonymous referees, the editor, and the participants at the Prospectus Workshop at the Department of Economics, Yale University, in 2002-2003 for helpful comments and suggestions. We also thank the participants at the Third International Conference on High Dimensional Probability, June 2002, and the 28th Conference on Stochastic Processes and Their Applications at the University of Melbourne, July 2002, where some of the results in the paper were presented.